\documentstyle{amsart}

\newcommand{\nc}{\newcommand}

\nc{\C}{{\cal C}}

\nc{\St}{\operatorname{St}^{\bullet}}
\nc{\B}{{\cal B}}
\nc{\A}{{\frak A}}

\nc{\N}{{\mathbf N}}
\nc{\CC}{{\mathbf C}}
\nc{\Z}{{\mathbf Z}}
\nc{\gr}{\operatorname{gr}}
\nc{\tr}{\operatorname{tr}}
\nc{\g}{{\frak g}}
\nc{\hatg}{\hat{\frak g}}
\renewcommand{\a}{{\frak a}}
\renewcommand{\b}{{\frak b}}
\nc{\gl}{{\frak g\frak l}}
\nc{\n}{{\frak n}}
\nc{\p}{{\frak p}}
\nc{\h}{{\frak h}}
\nc{\gothu}{{\frak u}}
\nc{\Ind}{\operatorname{Ind}}
\nc{\Coind}{\operatorname{Coind}}
\nc{\opp}{{\operatorname{opp}}}
\nc{\Ker}{\operatorname{Ker}}
\nc{\im}{\operatorname{Im}}
\nc{\Coker}{\operatorname{Coker}}
\nc{\dirlim}{\underset{\rightarrow}{\operatorname{lim}}}
\nc{\invlim}{\underset{\leftarrow}{\operatorname{lim}}}

\nc{\Ext}{\operatorname{Ext}^{\bullet}}
\nc{\ext}{\operatorname{Ext}}
\nc{\Tor}{\operatorname{Tor}_{\bullet}}
\nc{\tor}{\operatorname{Tor}}
\nc{\Tors}{\operatorname{Tor}_{\frac \infty 2+\bullet}}
\nc{\Exts}{\operatorname{Ext}^{\frac \infty 2+\bullet}}
\nc{\Hom}{\operatorname{Hom}^{\bullet}}
\nc{\ad}{\operatorname{ad}}
\renewcommand{\hom}{\operatorname{Hom}}

\renewcommand{\mod}{\operatorname{-mod}}
\nc{\Mod}{\operatorname{Mod}}

\renewcommand{\Bar}{\operatorname{Bar}}
\nc{\Barb}{\operatorname{Bar}^{\bullet}}

\nc{\upX}{X^{\uparrow}}
\nc{\upcD}{{\cal D}^{\uparrow}}
\nc{\upD}{D^{\uparrow}}
\nc{\dX}{X^{\downarrow}}
\nc{\dcD}{{\cal D}^{\downarrow}}
\nc{\dD}{D^{\downarrow}}
\nc{\upC}{{\cal C}^{\uparrow}}
\nc{\dC}{{\cal C}^{\downarrow}}
\nc{\underA}{\underline{A}}
\nc{\underB}{\underline{B}}
\nc{\underk}{\underline{k}}
\nc{\Db}{D^{\bullet}}
\nc{\map}{\longrightarrow}

\nc{\Q}{{\mathbf Q}}
\nc{\bs}{\bigskip\\}
\nc{\ms}{\smallskip\\}

\nc{\tilBar}{\widetilde{\operatorname{Bar}}}
\nc{\tilBarb}{\widetilde{\operatorname{Bar}}^{\bullet}}
\nc{\linBar}{\overline{\operatorname{Bar}}}
\nc{\linBarb}{\overline{\operatorname{Bar}}^{\bullet}}

\nc{\til}{\widetilde}
\nc{\oppA}{A^{\sharp}}

\nc{\Lemma}{{\bf Lemma:\ }}
\nc{\Theorem}{{\bf Theorem:}\ }
\nc{\Cor}{{\bf Corollary:}\ }
\nc{\Def}{{\bf Definition:}\ }
\nc{\Prop}{{\bf Proposition:}\ }
\nc{\Con}{{\bf Conjecture:}\ }
\nc{\Rem}{{\bf Remark:}\ }
\nc{\dok}{{\bf Proof.}\ }

\nc{\bul}{^{\bullet}}
\nc{\stand}{C^{\frac \infty 2+\bullet }}
\nc{\Cone}{\operatorname{Cone}}
\nc{\supp}{\operatorname{supp}}
\nc{\ten}{\otimes}

\nc{\ssn}{\subsection{}}
\nc{\sssn}{\subsubsection{}}

\nc{\hgt}{\operatorname{ht}}
\nc{\sqbinom}{\fracwithdelims[][0pt]}

\address{Independent University of Moscow, Pervomajskaya st. 16-18,
Moscow 105037, Russia}
\email{hippie@@ium.ips.ras.ru}

\thanks{Partially supported by Soros foundation.}

\author{S.~M.~Arkhipov}
\title{Semiinfinite cohomology of associative algebras and bar duality}
\begin{document}

\maketitle
\section{Introduction}
Semiinfinite cohomology of graded associative algebras was
first introduced in [Ar]. The setup for the definition
includes a graded associative algebra $A$ and two its graded
subalgebras $B$ and $N$ such that $N$ is positively graded,
$N$ is nonpositively graded and $A=B\ten N$ as a graded
vector space.

Basic examples  of this situation are provided by universal
enveloping algebras over graded Lie algebras. Let $\g$ be a
graded Lie algebra, and $A:=U(\g)$.  Then set
$\b:=\bigoplus_{n\le 0}\g_n$, $\n:=\bigoplus_{n>0}\g_n$, $B$
and $N$ are the universal enveloping algebras of the
corresponding Lie algebras.  The PBW Theorem provides the
{\em triangular decomposition}:  $U(\g)\cong U(\b)\ten
U(\n)$.

Lie algebra semiinfinite cohomology has a much longer
history. It appeared first in mathematics over 10 years ago
(see [F]). A proper homological construction for it  belongs
to Voronov (see [V]). Voronov treats Lie algebra
semiinfinite cohomology as an exotic derived functor of a
functor that is neither left nor right exact.

A similar approach for associative algebra semiinfinite
cohomology was presented in [Ar]. Let $k$ be a base field.
Consider a ``right semiregular $A$-module''
$S_A:=\Ind_N^A\Coind_k^N(k)$ and its endomorphism algebra
$\operatorname{End}_A(S_A)=:\oppA$. Suppose $\oppA$ is
augmented, i.~e.  there is a two-sided ideal
$\overline{A}^{\sharp}\in\oppA$ such that the quotient algebra
is $k$, let $\underk$ denote the trivial left
$\oppA$-module. Then semiinfinite cohomology with
coefficients in a graded $A$-module can be viewed as a ``two
sided derived functor'' of the functor
$\hom_{\oppA}(\underk,S_A\ten_A*)$ (see [V], [Ar]).

\ssn
This paper is devoted to an alternative construction of
associative algebra semiinfinite cohomology in terms of
graded associative algebras' bar duality.  Recall the idea
of the classical bar duality construction (see e.~g. [BGS]).
Let $A$ be a nonnegatively graded associative algebra such
that $\dim A_n<\infty$ and $A_0=k$. Then one can define a
canonical DG-algebra $A^{\vee}$ and a functor $D:\ A\mod\map
A^{\vee}\mod$ defined on suitably chosen derived categories
such that:  \begin{itemize} \item $D$ is an equivalence of
triangulated categories;
\item
$H\bul(A^{\vee})=\Ext_A(\underk,\underk)$; \item
$H\bul(A^{\vee\vee})=A$;
\item
$H\bul(D(M))=\Ext_A(\underk,M)$ for a graded $A$-module $M$;
\item
$H\bul(D^{-1}(N))=\Tor^{A^{\vee}}(\underk,N)$ for a $A^{\vee}$
DG-module $N$ satisfying certain grading conditions.
\end{itemize}
Recall also that $A^{\opp\vee}\cong A^{\vee\opp}$.

\ssn
The construction of associative algebra semiinfinite
cohomology is a direct generalization of the described one.

Suppose we have a graded associative algebra $A$ and its
nonpositively graded subalgebra $B$ satisfying certain
finiteness conditions. Suppose also that there is a left
$B$-augmentation of $A$, i.~e. a left $A$-module $\underB$
such that when restricted to $B$ the module $\underB$ is
isomorphic to the left regular representation of $B$. Note
that a triangular decomposition $A=B\ten N$ provides a
natural left $B$-augmentation of $A$:
$\underB:=A\ten_N\underk$.  In the second section we
construct a canonical DG-algebra $\Db(A,B)$ and  functors
$D^{\downarrow}_A:\ A\mod\map\Db(A,B)^\opp\mod$,
$D^{\uparrow}_A:\ A^\opp\mod\map\Db(A,B)\mod$ defined on
suitably chosen derived categories of modules such that
\begin{itemize}
\item
$D^{\downarrow}_A$ and $D^{\uparrow}_A$ are equivalences of
triangulated categories;
\item
$H\bul(\Db(A,B))=\Ext_A(\underB,\underB)$;
\item
$H\bul(D^{\uparrow}_A(M))=\Tor^A(\underB,M)$ for a right
$A$-module $M$.
\item
$H\bul(D^{\downarrow}_A(M))=\Ext_A(\underB,M)$ for a left
$A$-module $M$.
\end{itemize}
However it is no longer true
that taking the  opposite algebra commutes with duality. The
general case is as follows. There is another graded
associative algebra $\oppA$ containing $B^\opp$ as a
subalgebra (and having a triangular decomposition ---
$\oppA=B^\opp\ten N^\opp$ as a graded vector space --- if
$A$ has one) such that
$\Db(A,B)^\opp\cong\Db(\oppA,B^\opp)$.

In the third section we define semiinfinite $\ext$ functor
as follows. For a graded $A$-module $M$ and a graded
$\oppA$-module $L$ satisfying certain grading conditions
$$
  \Exts_A(L,M):=\Ext_{\Db(A,B)}(D^{\uparrow}_{\oppA}(L),D^{\downarrow}_A(M)).
$$
We prove also that the algebra $\oppA$ up to a completion coincides
with the algebra $\operatorname{End}_A(S_A)$ where $S_A:=\hom_B(A,B)$.
Note that in the case of triangular decomposition $A=B\ten N$ the $A$-module
$S_A$ is isomorphic to $\Ind_N^A\Coind_k^N\underk$. The algebra $\oppA$
plays the crucial part in all our considerations. In the Lie algebra case the
semiregular representation appeared first in the paper of Voronov [V],
the algebra $\oppA$ in the general case was introduced in [Ar]. We compare
the definition of semiinfinite $\ext$ functor with the one from
[Ar] and prove their coincidence.

The next  section is devoted to the semiinfinite cohomology in the
universal enveloping algebra case. Suppose we have a graded Lie
algebra $\g$ and its triangular decomposition $\g=\b\oplus\n$ as a
graded vector space. Consider the universal enveloping algebra
$A:=U(\g)$ with the induced triangular deccomposition. We prove that
the algebra $\oppA$ is also a universal enveloping algebra of a Lie
algebra that differs from $\g$ by a one-dimensional central extension.
This central extension is defined with the help of a 2-cocycle of $\g$
known as the {\em critical cocycle} (see e.~g.~[FFr]). Thus we prove
that in the Lie algebra case our semiinfinite cohomology of the
universal enveloping algebra coincides with Lie algebra semiinfinite
cohomology (see e.~g~[F]).  We conclude the section with a direct
calculation of the critical 2-cocycle for affine Lie algebras and for
different triangular decompositions. It turns out that the cocycle
itself depends on the type of the decomposition, still cocycles for
different decompositions define the same cohomology class.

The author would like to thank M.~Finkelberg and L.~Positselsky for
helpful discussions.  L.~Positselsky also showed  the unpublished
paper [P] to the author, it stimulated the whole investigation a lot.
B.~Feigin introduced the author to the subject and explained both
mathematical and physical pictures of it with an unbelievable
patience.  This has been a matter of great importance for me and my
gratitude cannot even be expressed properly.

\subsection{Conventions}
In most part of the paper we work over an arbitrary field
$k$. All the vector spaces are considered over that field.
In the end of the fourth section our base field is $\mathbf{C}$.

\section{Relative bar construction}
\ssn
\label{setup}%
Suppose we have a graded associative
algebra $A=\underset{n\in \Z}{\bigoplus}A_n$. Let $B$ and $N$ be
graded subalgebras
in $A$ satisfying
the following conditions:

 \qquad(i) $N$ is positively graded;

 \qquad(ii) $N_0=k;$

 \qquad(iii) $\dim N_n<\infty$ for any $n\in\N;$

 \qquad(iv) $B$ is negatively graded;

 \qquad(v) the multiplication in $A$ defines the isomorphisms of graded
 vector spaces
$$
 B\otimes N\map A\text{ and }N\otimes B\map A.
$$
In particular $N$ is naturally augmented.
We denote the augmentation ideal
$\underset{n>0}{\bigoplus}N_n$ by $\overline{N}$.

The category of left graded $A$-modules with morphisms that
preserve gradings is denoted by $A\mod$. We define the
functor of the grading shift
$$
 A\mod\map A\mod:\ M\mapsto
 M\langle i\rangle,\ M\langle i\rangle_m:=M_{i+m}, i\in\Z.
$$
The space $\underset
{i\in\Z}{\bigoplus}\hom_{A\mod}(M_1,M_2\langle i\rangle)$ is
denoted by $\hom_A(M_1,M_2)$.

\ssn
We fix a left $B$-augmentation on $A$ provided by the
isomorphism of graded left $B$-modules $B\cong
A\ten_N\underk$ where $\underk:=N/\overline{N}$ is the
trivial $N$-module. The $A$-module $A\ten_Nk$ is denoted by
$\underB$.

\ssn
We introduce certain subcategories in the category of
complexes ${\cal K}om(A\mod)$.  For
$M\bul\in{\cal K}om(A\mod)$ the support of $M$ is defined as
follows:
$$
 \supp M\bul:=\{(p,q)\in\Z^{\oplus 2}|M_p^q\ne
 0\}.
$$
For $s_1,s_2,t_1,t_2\in \Z,\ s_1,s_2>0$, the set $
\{(p,q)\in \Z^{\oplus 2}|s_1q+p\ge t_1,\ s_2q-p\ge t_2\} $
(resp. the set $ \{(p,q)\in \Z^{\oplus 2}|s_1q+p\le t_1,\
s_2q-p\le t_2\}) $ is denoted by $\upX(s_1,s_2,t_1,t_2)$
(resp. by $\dX(s_1,s_2,t_1,t_2)$).

\sssn
Let $\upC(A)$ (resp. $\dC(A)$) be the full subcategory in
${\cal K}om(A\mod)$ consisting of complexes $M\bul$ that
satisfy the following condition:

\quad(U) there exist $s_1,s_2,t_1,t_2\in \Z,\ s_1,s_2>0$,
such that $\supp M\bul\subset\upX(s_1,s_2,t_1,t_2)$ (resp.

\quad(D) there exist $s_1,s_2,t_1,t_2\in \Z,\ s_1,s_2>0$,
such that $\supp M\bul\subset\dX(s_1,s_2,t_1,t_2)$).

\sssn
Recall the construction of standard relative bar resolution
of a complex of $A$-modules $M\bul$.

The standard bar resolution $\tilBarb (A,B,M)\in {\cal
K}om(A\mod)$ of an $A$-module $M$ is defined as follows:
\begin{gather*}
  \tilBar^{-n}(A,B,M):=A\otimes_B\ldots \otimes_BA\otimes_BM \
  (n+1 \text{ times}), \\
  \operatorname{d}(a_0\otimes \ldots \otimes a_n\otimes v)=
  \sum_{s=0}^{n-1}(-1)^s a_0\otimes\ldots\otimes
  a_sa_{s+1}\otimes\ldots\otimes v
  +(-1)^na_0\otimes\ldots\otimes a_{n-1}\otimes a_nv.
\end{gather*}
Here
$a_0,\ldots ,a_n\in A,\ v\in M$.

\sssn
\Lemma The subspace $\linBarb (A,B,M):$
$$
  (\linBar)^{-n}(A,B,M):=\{ a_0\otimes\ldots\otimes a_n\otimes
  v\in \tilBar^{-n}(A,B,M)|\ \exists \ s\in\{1,\ldots,n\}:\
  a_s\in B\}
$$
is a subcomplex in $\tilBarb (A,B,M)$.\qed
\ms
The quotient $\Barb (A,B,M):=\tilBarb (A,B,M)/\linBarb
(A,B,M)$ is called the restricted bar resolution of the
$A$-module $M$ with respect to the subalgebra $B$.

For a complex of $A$-modules $M\bul\in\upC(A)$ its relative
restricted bar resolution is defined as a total complex of
the bicomplex $\Barb (A,B,M\bul)$.

\sssn
\Lemma Let $M\in \upC(A)$. Then

\qquad(i) $\Bar^{-m}(A,B,M\bul)=N\ten\overline{N}^{\ten
m}\ten M\bul$ as a vector space;

\qquad(ii) $\Barb (A,B,M\bul)\in\upC(A)$;

\qquad(iii) $\Barb (A,B,M\bul)$ is quasiisomorphic to $M$ as
a complex of $A$-modules.  \qed

Note that for $M\in A\mod$ the complex $\Barb (A,B,M)$ is
isomorphic to $\Barb(N,k,M)$ as a complex of $N$-modules.
Here $k$ is considered as a subalgebra in $N$.

\ssn
Consider the relative restricted bar resolution $\Barb
(A,B,\underB)$ of the $A$-module $\underB$ and the complex
of graded vector spaces
$$
 \Db(A,B):=\Hom_A(\Barb(A,B,\underB),\underB).
$$
Clearly $\Db(A,B)\cong
\Hom_k(\underset{n>0}{\bigoplus}\overline{N}^{\ten
n},\underB)$ as a vector space, and $D(A,B)$ belongs to
$\dC({\cal V}ect)$.

\sssn \label{mult}%
We introduce a structure of a DG-algebra on $\Db(A,B)$.
First we define a DG-algebra structure on
$\Hom_A(\tilBarb(A,B,\underB),\underB)$. Note that by
Schapiro lemma
$$
  \Hom_A(\tilBar^{-m}(A,B,\underB),\underB)=
  \Hom_B(A\ten_B\ldots\ten_BA\ten_B\underB,\underB)\ (m\text{ times}).
$$
Let $f\in \hom^m_A(\tilBarb(A,B,\underB),\underB),\ g\in
\hom^n_A(\tilBarb(A,B,\underB),\underB)$, i.~e.
$$
  f:\
  \underbrace{A\ten_B\ldots\ten_BA}_m\ten_B\underB\map
  \underB,\ g:\
  \underbrace{A\ten_B\ldots\ten_BA}_n\ten_B\underB\map
  \underB,
$$
both $f$ and $g$ are $B$-linear. By definition
set $f\cdot g:\
\underbrace{A\ten_B\ldots\ten_BA}_{m+n}\ten_B\underB\map
  \underB:$
$$
  (f\cdot g)(a_1\ten\ldots\ten a_{m+n}\ten b):=
  f(a_1\ten\ldots\ten a_n\ten g(a_{n+1}\ten\ldots\ten
  a_{m+n}\ten b)).
$$
\sssn
\Lemma

\qquad(i) The multiplication
equips $\Hom_A(\tilBarb(A,B,\underB),\underB)$ with a
DG-algebra structure;
$$
 (\text{ii})\qquad{}\ \
  \Db(A,B)=\{f\in\Hom_A(\tilBarb(A,B,\underB),\underB)\ |
  f\equiv 0 \text{ on }
  \linBarb(A,B,\underB)\subset\tilBarb(A,B,\underB)\}
$$
is a  DG-subalgebra in
$\Hom_A(\tilBarb(A,B,\underB),\underB)$.\qed

\sssn
\Rem
(i)
$D^0(A,B)=\hom_A(A\ten_B\underB,\underB)=
\hom_B(\underB,\underB)\cong B^{\opp}$ as an algebra (yet
the inclusion $B^{\opp}\hookrightarrow \Db(A,B)$ is not a
morphism of DG-algebras --- the differential in $\Db(A,B)$
does not preserve $B^{\opp}$).

(ii) Consider the induction functor $\Ind_N^A:\ N\mod\map
A\mod$. Then the canonical map
\begin{multline*}
  \Hom_N(\Barb(N,k,\underk),\underk)\map\Hom_A(\Ind_N^A(\Barb(N,k,\underk)),
  \Ind_N^A(\underk))\\
  \cong\Hom_A(\Barb(A,B,\underB),\underB)=\Db(A,B)
\end{multline*}
is an inclusion of DG-algebras.

(iii) Denote its image by $\Db(N,k)$. Then
$
\Db(A,B)=\Db(N,k)\ten B^{\opp}
$
as a graded vector space.

\sssn \label{functor}%
Let  $M\bul\in\dC(A),\ M^{\prime\bullet}\in\upC(A^\opp)$.
By definition set
$$
  \dD(A,B,M\bul):=\Hom_A(\Barb(A,B,\underB),M\bul),\
  \upD(A,B,M^{\prime\bullet}):=M^{\prime\bullet}\ten_A\Barb(A,B,\underB).
$$
Evidently the vector space $\dD(A,B,M\bul)$
 (resp. $\upD(A,B,M^{\prime\bullet})$) belongs to
$\dC({\cal V}ect)$ (resp. to $\upC({\cal V}ect)$).  Similarly
to~\ref{mult} we define the right action of $\Db(A,B)$ on
 $\dD(A,B,M\bul)$ (resp the left action of $\Db(A,B)$ on
$\upD(A,B,M^{\prime\bullet})$) as follows.

Recall that
$\Hom_A(\tilBar^{-m}(A,B,\underB),M\bul)\cong
\Hom_B(\underbrace{A\ten_B\ldots\ten_BA}_m\ten_B\underB,M\bul).$

Suppose we have
$f:\ \underbrace{A\ten_B\ldots\ten_BA}_m\ten_B\underB\map M\bul,\
g:\ \underbrace{A\ten_B\ldots\ten_BA}_n\ten_B\underB\map\underB,$
both $f$ and $g$ are $B$-linear. Then by definition set
$(f\cdot g):\
\underbrace{A\ten_B\ldots\ten_BA}_{m+n}\ten_B\underB\map
 M\bul:$
$$
  (f\cdot g)(a_1\ten\ldots\ten a_{m+n}\ten b):=
  f(a_1\ten\ldots\ten a_m\ten g(a_{m+1}\ten\ldots\ten
  a_{m+n}\ten b)).
$$
Note also that
$M^{\prime\bullet}\ten_A\tilBarb(A,B,\underB)\cong
 \bigoplus_{n\le 0}M'\ten_B(
\underbrace{A\ten_B\ldots\ten_BA}_{-n})\ten_B\underB.$

Suppose we have
$$
 x=m\ten a_1\ldots\ten a_n\ten b\in M^{\prime\bullet}\ten_B
 \underbrace{A\ten_B\ldots\ten_BA}_n\ten_B\underB,\ f:\
 \underbrace{A\ten_B\ldots\ten_BA}_m\ten_B\underB\map\underB,
$$
and $f$ is $B$-linear. Then by definition set
$f\cdot x\in M^{\prime\bullet}\ten_B(
\underbrace{A\ten_B\ldots\ten_BA}_{n-m})\ten_B\underB:$
$$
 f\cdot x:=m\ten a_1\ten\ldots\ten a_{n-m}\ten f(a_{n-m+1}
 \ten\ldots\ten a_n\ten b).
$$
\sssn
\Lemma

\qquad(i) The multiplication equips
$\Hom_A(\tilBarb(A,B,\underB),M\bul)$ with a structure of
the right DG-module

\qquad over $\Db(A,B)$;

\qquad(ii)
$\dD(A,B,M\bul)\subset\Hom_A(\tilBarb(A,B,\underB),M\bul)$
is a DG-submodule;

\qquad(iii) the multiplication equips
$M^{\prime\bullet}\ten_A\tilBarb(A,B,\underB)$ with a
structure of the left DG-module over

\qquad $\Db(A,B)$;

\qquad(iv) $\upD(A,B,M^{\prime\bullet})$ is a quotient
DG-module of $M^{\prime\bullet}\ten_A\tilBarb(A,B,\underB)$.
\qed

\sssn
Denote the category of right (resp. left) DG-modules
$$
  M\bul=\underset{p,q\in\Z}{\bigoplus}M_p^q,\
  \operatorname{d}:\ M_p^q\map M_p^{q+1},
$$
over $\Db(A,B)$, with morphisms being
morphisms of DG-modules that preserve also the second
grading, by $\Db(A,B)\mod$ (resp. by $\Db(A,B)^\opp\mod$).
The subcategory in $\Db(A,B)\mod$ (resp. in
$\Db(A,B)^\opp\mod$) that consists of DG-modules satisfying
the condition (D) (resp. (U)) is denoted by $\dC(\Db(A,B))$
(resp.  by $\upC(\Db(A,B)^\opp)$.

The localizations of $\dC(A)$, $\upC(A^\opp)$,
$\dC(\Db(A,B))$, $\dC(\Db(A,B))$
and $\upC(\Db(A,B)^\opp)$ by the class of quasiisimorphisms
are denoted by $\dcD(A)$, $\upcD(A^\opp)$, $\dcD(\Db(A,B))$,
and $\upcD(\Db(A,B)^\opp)$ respectively.

By~\ref{functor} $\dD$ and $\upD$ define the  functors
$$
 \dD_A:\ \dC(A)\map \dC(\Db(A,B))\text{ and }\upD_A:\
 \upC(A^\opp)\map \upC(\Db(A,B)^\opp).
$$
\sssn
\Theorem \label{equiv}%

\qquad(i) The functor $\dD_A$ is well defined as a functor
from $\dcD(A)$ to $\dcD(\Db(A,B))$;

\qquad(ii) $\dD_A:\ \dcD(A)\map \dcD(\Db(A,B))$ is an
equivalence of triangulated categories;

\qquad(iii) the functor $\upD_A$ is well defined as a
functor from $\upcD(A^\opp)$ to $\upcD(\Db(A,B)^\opp)$;

\qquad(iv) $\upD_A:\ \upcD(A^\opp)\map \upcD(\Db(A,B)^\opp)$
is an equivalence of triangulated categories.

\subsection{Proof of Theorem~\ref{equiv}}

We present here the proof of the first two statements, the
other two are verified in a similar way.

First we prove that $\dD_A$ maps quasiisomorphisms into
quasiisomorphisms.

\sssn
\Lemma \label{deriv}%
Let $M\bul\in\dC(A)$ be exact. Then
$\dD_A(M\bul)\in\dC(\Db(A,B))$ is also exact.

\dok
As a complex of vector spaces
$\dD_A(M\bul)=\Hom_N(\Barb(N,k,\underk),M\bul)$.  Consider
the spectral sequence of this bicomplex.
$$
  E_0^{p,q}=\hom_k((\overline{N})^{\ten -p},M^q)=
  \underset{n\in\Z}{\bigoplus}\,\underset{\sum
  t_j-s=n}{\bigoplus} \hom_k(N_{t_1}\ten\ldots\ten
  N_{t_{-p}},M_s^q)
$$
As $\dD_A(M\bul)$ satisfies  the
condition (D), in a fixed grading the spectral sequence is
situated in an area of the $(p,q)$ plane bounded both in $p$
and $q$.  Thus the spectral sequence converges. But the
first term of it looks as follows.
$$
  E_1^{p,q}=\hom_k((\overline{N})^{\ten -p},H^q(M)).
$$
In particular $E_1^{p,q}=0$ if $M\bul$ is exact.\qed

The first statement of the theorem is proved. To prove the
second one we will construct a functor
$\overline{D}^{\downarrow}_A:\ \dC(\Db(A,B))\map \dC(A)$.

Consider the complex
$K\bul:=\Hom_A(\Barb(A,B,\underA),\underB)$. As a vector
space
$$
 K\bul=\Hom_k(\underset{n\in\Z}{\bigoplus}N\ten
 \overline{N}^{\ten n},\underB).
$$
In particular
$K\bul\in\dC({\cal V}ect)$. We define the {\em left\/} action
of $\Db(A,B)$ on $K\bul$ similarly to~\ref{mult}.  The left
action of $A$ on $K\bul$ is defined using the right
multiplication in the $A$-$A$ bimodule $\underA$.

\sssn
\Lemma
$K\bul$ is a $\Db(A,B)$-$A$ DG-bimodule.\qed

Note that as a $\Db(A,B)$-module
$K\bul\cong\Db(A,B)\ten_NN^*$ (the differentials both in
$\Db(A,B)$ and $K\bul$ are forgotten).

Consider a functor
$
  \overline{D}^{\downarrow}_A:\ \Db(A,B)\mod\map A\mod,\
  X\bul\mapsto K\bul\ten_{\Db(A,B)}X\bul.$

\sssn
\Lemma

\qquad(i) $\overline{D}^{\downarrow}_A:\
\dC(\Db(A,B)\map\dC(A)$;

\qquad(ii) $\overline{D}^\downarrow_A$ is well defined on
the corresponding derived categories.
\ms
\dok
(i) Note that $\overline{D}^{\downarrow}_A(X\bul)=N^*\ten X\bul$ as a
vector space.

(ii) The proof is parallel to the one of~\ref{deriv}.\qed

Note that for $M\in A\mod$ the complex
$\Hom_A(\Barb(A,B,\underA),M)$ is a complex of left
$A$-modules. The multiplication by elements of $A$ is
provided by the right multiplication in the $A$-bimodule
$\Barb(A,B,\underA)$.

Note also that the complex
$\Hom_A(\Barb(A,B,\Barb(A,B,\underB)),\underB)$ is a
DG-bimodule over $\Db(A,B)$. The construction of the right
and the left action of $\Db(A,B)$ on the complex is as
follows.
\begin{multline*}\label{doublebar}%
 \Bar^{-p}(A,B,\Bar^{-q}(A,B,\underB))=\\
 A\ten_B\left(\underbrace{(A/B)\ten_B\ldots\ten_B(A/B)}_p\right)
 \ten_B A\ten_B\left(\underbrace{(A/B)
 \ten_B\ldots\ten_B(A/B)}_q\right)\ten_B\underB.
\end{multline*}
Suppose we have
\begin{gather*}
  f_1\in\Hom_A(\Bar^{-s}(A,B,\underB),\underB)=
  \Hom_B(\underbrace{A\ten_B\ldots\ten_BA}_{s}\ten_B\underB,\underB),\\
  f_2\in\Hom_A(\Bar^{-t}(A,B,\underB),\underB)=
  \Hom_B(\underbrace{A\ten_B\ldots\ten_BA}_{t}\ten_B\underB,\underB),\\
  g\in\Hom_A(\Bar^{-p}(A,B,\Bar^{-q}(A,B,\underB),\underB)\hookrightarrow
  \Hom_B(\underbrace{A\ten_B\ldots\ten_BA}_{p+q+1}\ten_B\underB,\underB).
\end{gather*}
Then set
\begin{multline*}
  (f_1\cdot g\cdot
  f_2)(a_1\ten\ldots\ten a_{s+t+p+q+1}\ten b)=\\
  f_1(a_1\ten\ldots\ten a_{s}\ten g(a_{s+1}\ten\ldots\ten
  a_{s+p+q+1}\ten f_2(a_{s+p+q+2}\ten\ldots\ten
  a_{s+t+p+q+1}\ten b))).
\end{multline*}
One can check directly
that the differential in the complex satisfies the Leibnitz
rule both for the left and  the right DG-module structures.

Consider a subspace
$L\bul\subset\Barb(A,B,\Barb(A,B,\underB))$ defined as
follows.
\begin{gather*}
  L^{-p,-q}:=A\ten_B\left(\underbrace{(A/B)\ten_B\ldots
  \ten_B(A/B)}_p\right)
  \ten_BB\ten_B \left(\underbrace{(A/B)
  \ten_B\ldots\ten_B(A/B)}_q\right)\ten_B\underB,\\
  L\bul:=\underset{p,q\in\Z}{\bigoplus}L^{p,q}\hookrightarrow
  \Barb(A,B,\Barb(A,B,
  \underB)).
\end{gather*}
\sssn
\Lemma
$L\bul$ is a subcomplex in $\Barb(A,B,\Barb(A,B,\underB))$.  \qed

\sssn
\Lemma

\qquad(i) Let $M\bul\in\dC(A)$. Then
$\overline{D}^{\downarrow}_A\circ
D^{\downarrow}_A(M\bul)\cong\Hom_A(\Barb(A,B,\underA),M\bul)$ as a
complex of

\qquad $A$-modules.

\qquad(ii) Let $X\bul\in\dC(\Db(A,B))$. Then
$$
  D^{\downarrow}_A\circ\overline{D}^{\downarrow}_A(X\bul)\cong
  X\bul\ten_{\Db(A,B)}\Hom_A(\Barb(A,B,\Barb(A,B,\underB)),\underB).\qed
$$
The canonical morphism of $A$-modules
$\Barb(A,B,\underA)\map \underA$ provides a morphism of
functors $\phi:\
Id_{\dC(A)}\map\overline{D}^{\downarrow}_A\circ
D^{\downarrow}_A$.

The canonical morphism of complexes
$\Barb(A,B,\underB)\map\Barb(A,B,\Barb(A,B,\underB))$:
 \begin{gather*}
 \Bar^{-m}(A,B,\underB)\map\underset{p+q=m}{\bigoplus}L^{-p,-q}
 \subset\Barb(A,B,\Barb(A,B,\underB))\\
 a_1\ten\ldots\ten a_{m+1}\ten b\mapsto \sum_{i=1}^{m+1}
 a_1\ten\ldots\ten a_i\ten 1\ten a_{i+1}\ten\ldots \ten
 a_{m+1}\ten b \end{gather*} provides a morphism of
DG-bimodules over $\Db(A,B)$
$$
  \Hom_A(\Barb(A,B,\Barb(A,B,\underB)),\underB)\map\Db(A,B)
$$
that leads to the morphism of functors $\psi:\
D^{\downarrow}_A\circ\overline{D}^{\downarrow}_A
\map\operatorname{Id}_{\dC(\Db(A,B))}$.

\sssn
\Lemma
Both $\phi$ and $\psi$ are quasiisomorphisms.
\qed

Thus Theorem~\ref{equiv} is proved.\qed

\sssn
\Rem
To construct an  inverse functor for $\upD_A$ one has to
consider a functor
$$
  \overline{D}^\uparrow_A:\
  \upC(\Db(A,B)^\opp)\map\upC(A^\opp),\
  \overline{D}^\uparrow_A:\ X\bul\mapsto
  \Hom_{\Db(A,B)^\opp}(K\bul,X\bul).
$$

\section{The algebra $\oppA$ and semiinfinite cohomology}

In this section we give a realization of semiinifinite
cohomology in terms of bar duality. The classical
construction of that type was obtained in [BGG].  The Koszul
duality functor $K$ constructed there gave an equivalence of
suitably chosen derived categories of modules over a
quadratic Koszul algebra $A$ and its quadratic dual algebra
$A^!$ (see [BGG], [BGS]). The counterpart of the classical
$\ext$ functor in that construction was given by the
following statement.

\ssn
\Prop
Let $M$ be a graded module over a quadratic Koszul algebra
$A$. Then
$H\bul(K(M))=\underset{n\in\Z}{\bigoplus}
\Ext_{A\mod}(\underk,M\langle n\rangle)$. Here as before
$\underk$ denotes the trivial module over $A$ placed in the
zero grading and $\langle*\rangle$ denotes the shift in the
category of graded modules.\qed

For a graded algebra $A$ and  graded $A$-modules $M_1$ and
$M_2$ we denote the space
$\underset{n\in\Z}{\bigoplus}\Ext_{A\mod}(M_1,M_2\langle
n\rangle)$ by $\Ext_A(M_1,M_2)$.

\ssn
We return to the situation of the previous section.

Let $\til{\cal C}^{\downarrow}(\Db(A,B))$ be the category of
DG-modules over $\Db(A,B)$ whose cohomology satisfies the
condition (D), $\til{\cal D}^{\downarrow}(\Db(A,B))$ denotes
its localization by the class of quasiisomorphisms.

\sssn
\Lemma            \label{truncation}%
The natural inclusion functor provides the equivalence of
categories
$$
  \dcD(\Db(A,B))\til{\map}\til{\cal
  D}^\downarrow(\Db(A,B)).
$$
\ms
\dok
We construct
truncation functors similar to those described in [GM] for
the usual categories of complexes over abelian categories.

First note that for any $M\bul\in\Db(A,B)\mod$ and an
arbitrary integer $n$ the complex
$$
 \sigma_{\le n}M\bul:\
  \sigma_{\le n}M_p^q:= \begin{cases} M_p^q, &\text{if } p\le
    n,\\ 0, &\text{if } p>n \end{cases}
$$
is a DG-submodule in $M\bul$.

Denote the set
$
  \{(p,q)\in\Z^{\oplus 2}\setminus\dX(s_1,s_2,t_1,t_2)|
  (p+1,q)\in\dX(s_1,s_2,t_1,t_2)\}.
$
(resp. the set
$
  \{(p,q)\in\Z^{\oplus 2}\setminus\dX(s_1,s_2,t_1,t_2)|
  (p-1,q)\in\dX(s_1,s_2,t_1,t_2)\}).
$
by $X^{\swarrow}(s_1,s_2,t_1,t_2)$
(resp. by $X^{\searrow}(s_1,s_2,t_1,t_2)$).

Choose $s_1,s_2,t_1,t_2$ such that
$$
  \supp H\bul(M\bul)\subset\dX(s_1,s_2,t_1,t_2)
  \text{ and }\supp\Db(A,B)\subset\dX(s_1,s_2,0,0).
$$
Let $(p_0,q_0)$ be the vertex of the convex cone
$\dX(s_1,s_2,t_1,t_2)$.  Consider a submodule
$\tau_{\swarrow}(s_1,s_2,t_1,t_2)(M\bul)$ in $\sigma_{\le
p_0}M\bul$:
$$
  \tau_{\swarrow}(s_1,s_2,t_1,t_2)(M\bul)_p^q:=
  \begin{cases} 0, &\text{if }p>p_0\\
     \Ker(\operatorname{d}_p^q), &\text{if }(p,q)\in
     X^{\swarrow}(s_1,s_2,t_1,t_2)\\ M_p^q,
  &\text{otherwise.} \end{cases}
$$
One can check directly
that the defined submodule is quasiisomorphic to $M\bul$.

Consider also a  submodule
$\tau_{\searrow}(s_1,s_2,t_1,t_2)(M\bul)$ in $\sigma_{\le
p_0}M\bul$:
$$
 \tau_{\searrow}(s_1,s_2,t_1,t_2)(M\bul)_p^q:=
  \begin{cases} 0, &\text{if }p>p_0\\
     \im(\operatorname{d}_p^q), &\text{if }(p,q)\in
     X^{\searrow}(s_1,s_2,t_1,t_2)\\ M_p^q,
  &\text{otherwise.} \end{cases}
$$ One can check directly
that the defined submodule is quasiisomorphic to $M\bul$.

But
$\tau_{\searrow}(s_1,s_2,t_1,t_2)
\tau_{\swarrow}(s_1,s_2,t_1,t_2)M\bul$ satisfies the
condition (D).  Thus every DG-module from $\til{\cal
C}^{\downarrow}(\Db(A,B))$ is quasiisomorphic to a DG-module
from $\dC(\Db(A,B))$.

That means that the inclusion functor
$\dcD(\Db(A,B))\map\til{\cal D}^\downarrow(\Db(A,B))$ is
surjective on classes of isomorphism of objects. On the
other hand let
$$
 M_1\bul,M_3\bul\in\dC(\Db(A,B)),\
  M_2\bul\in\til{\cal C}^{\downarrow}(\Db(A,B)),\
  M_1\bul\overset{s}{\longleftarrow}M_2\bul\overset{f}{\map}M_3\bul,
$$
where $s$ is a quasiisomorphism, be a diagram in $\til{\cal
C}^{\downarrow}(\Db(A,B))$ representing a morphism in
$\til{\cal D}^{\downarrow}(\Db(A,B))$ (see [GM] for the
explicit description of morphisms in  derived categories).
Then there exist $s_1,s_2,t_1,t_2$ such that the diagram is
equivalent to the following one:
$$
 M_1\bul\longleftarrow
 \tau_{\swarrow}(s_1,s_2,t_1,t_2)\tau_{\searrow}(s_1,s_2,t_1,t_2)M_2\bul
 \map M_3\bul.
$$
That means that the inclusion functor on
derived categories is surjective on the class of morphisms.
In the same way one checks that if two diagrams represent
the same morphism in $\til{\cal D}^{\downarrow}(\Db(A,B))$
then the corresponding truncated diagrams are also
equivalent in $\dD(\Db(A,B))$. Thus the inclusion functor is
an equivalence of categories.  \qed

\sssn
\Lemma \label{homotopy}%
Let $P^{\bullet,\bullet}$ be a  bicomplex over an abelian
category $\cal A$ such that for $r>>0$ all $P^{r,\bullet}=0$
and for every $r\in\Z$ the complex $P^{r,\bullet}$ is
homotopic to zero.  Then the total complex $P\bul $ is also
homotopic to zero.
\ms
\dok
The differential in the
bicomplex $P^{\bullet,\bullet}$ is the sum of two components
of the bigradings $(0,1)$ and $(1,0)$ respectively:
$d=d_1+d_2$. Choose a homotopy map
$$
 h_r:\ P^{r,\bullet}\map P^{r,\bullet-1};\
 d_1h_r+h_rd_1=\operatorname{Id}_{P^{r,\bullet}}.
$$
Set $h:=\sum h_r$. Then the map
$dh+hd=\operatorname{Id}_{P^{\bullet\bullet}}+d_2h+hd_2$ is
the sum of an invertible map and a nilpotent one. Thus it is
invertible itself.  \qed

We say that two morphisms of DG-modules over $\Db(A,B)$ are
homotopic to each other:  $f_1\sim f_2,\ f_1,f_2:\
M_{(1)}\bul\map M_{(2)}\bul$, if there exists a morphism of
bigraded $\Db(A,B)$-modules
$$
 h:\ M_{(1)}\bul\map
 M_{(2)}\bul,\ h_p^q:\ M_{(1)p}^q\map M_{(2)p}^{q-1},
$$ such
that $f_1-f_2=dh+hd$. Two DG-modules $M_1$ and $M_2$ are
homotopically equivalent, $M_1\sim M_2$, if there exist
morphisms of DG-modules
$$
 f:\ M_1\map M_2 \text{ and } g:\
 M_2\map M_1:\ fg\sim \operatorname{Id}_{M_2},\
 gf\sim\operatorname{Id}_{M_1}.
$$
A DG-module $M\bul$ is homotopic to zero if
$\operatorname{Id}_M\sim 0$.

Now we are to  construct a realization of the category
$\dcD(\Db(A,B))$ similar to the realization of the usual
bounded derived category of an abelian category as the
homotopical category of complexes consisting of projective
objects.

\sssn
\Lemma \label{bar0}%
Let $M\bul\in\Db(A,B)\mod$ be an exact DG-module. Then the
DG-module \\$\Barb(\Db(A,B),k,M\bul)$ is homotopic to zero.
\ms
\dok
Consider a bicomplex $\Barb(\Db(A,B),k,M\bul)$:
$$
  \Bar^{u,v}(\Db(A,B),k,M\bul):=\underset{r}{\bigoplus}\underset
  {\underset{p_i\ne 0,p+\sum p_i+s=r} {u=-n,q+\sum
  q_i+t=v}}{\bigoplus}D_p^q(A,B)\ten D_{p_1}^{q_1}(A,B)
  \ten\ldots\ten D_{p_n}^{q_n}(A,B)\ten M_s^t.
$$
The differential in the bicomplex is a sum of two components
$d_1$ and $d_2$.  The bicomplex \\
$\Bar^{u,v}(\Db(A,B),k,M\bul)$ is bounded from above in $u$.

Note that each of the complexes of vector spaces
$\overline{D}\bul(A,B)^{\ten p}\ten M\bul$ is exact.
Choose  homotopy maps $h_p$ for these complexes.
Then $\operatorname{Id}\ten h_u$ will be
a homotopy map in the  $\Db(A,B)$-module
$\Db(A,B)\ten \overline{D}(A,B)^{\ten u}\ten M\bul$
with the component of the differential $d_2$ forgotten. By
Lemma~\ref{homotopy} the  DG-module
$\Barb(\Db(A,B),k,M\bul)$ is itself homotopic to zero.  \qed

\sssn
\Cor \label{bar1}%
For $M\bul\in\dC(\Db(A,B))$
$$
 \Barb(\Db(A,B),k,M\bul)\sim
 \Barb(\Db(A,B),k,\Barb(\Db(A,B),k,M\bul)).\qed
$$
The
homotopical category ${\cal K}(\dC(\Db(A,B)))$ is the
category with the class of objects the same as in
$\dC(\Db(A,B))$ and morphisms being classes of homotopical
equivalence of morphisms in $\dC(\Db(A,B))$. Standard
considerations (see [GM]) show that it is a triangulated
category.

Consider the full subcategory ${\cal BAR}(\Db(A,B))$ in the
category ${\cal K}(\Db(A,B)\mod)$, a DG-module
$M\bul\in{\cal BAR}(\Db(A,B))$ iff $M\bul$ is isomorphic to
some DG-module of the type $\Barb(\Db(A,B),k,M')$. Then
Corollary~\ref{bar1} implies that ${\cal BAR}(\Db(A,B))$ is a
triangulated category: one has to check that for every
morphism of DG-modules
$$
  f:\ \Barb(\Db(A,B),k,M\bul)\map
  \Barb(\Db(A,B),k,M^{\prime\bullet})
$$
the DG-module $\Cone f$
is also homotopically equivalent to a DG-module in
${\cal BAR}(\Db(A,B))$. Now consider a morphism
$$
  \Bar(f):\ \Barb(\Db(A,B),k,\Barb(\Db(A,B),k,M\bul))
  \map \Barb(\Db(A,B),k,\Barb(\Db(A,B),k,M^{\prime\bullet}))
$$
that corresponds to $f$. Since ${\cal K}(\Db(A,B)\mod)$ is a
triangulated category itself, $\Cone(f)\sim\Cone(\Bar(f))$.
But $\Cone(\Bar(f))\cong\Barb(\Db(A,B),k,\Cone(f))$ and we
are done.

Let ${\cal D}(\Db(A,B))$ be the localization of $\Db(A,B)\mod$
by the class of quasiisomorphisms,

\sssn
\Prop
The functor of localization by the class of
quasiisomorphisms provides an equivalence of triangulated
categories $L:\
{\cal BAR}(\Db(A,B))\widetilde{\map}{\cal D}(\Db(A,B))$.

\dok
Let us construct the inverse functor. Consider a functor
$$
  \overline{L}:\ \Db(A,B)\mod\map{\cal BAR}(\Db(A,B)),\
  M\bul\mapsto \Barb(\Db(A,B),k,M\bul).
$$
By
Lemma~\ref{bar0} the functor $\overline{L}$ takes exact
DG-modules to zero, thus quasiisomorphisms are mapped into
quasiisomorphisms. By the characteristic property of the
derived category $\overline{L}$ induces a functor
${\cal D}(\Db(A,B))\map {\cal BAR}(\Db(A,B))$.  This functor
is the one we need.  \qed

\sssn
\Lemma  \label{HHH}%
Let $M\bul,M^{\prime\bullet}\in\Db(A,B)\mod$. Then
$$
  \hom_{{\cal K}(\Db(A,B))}(\Barb(\Db(A,B),k,M\bul),M^{\prime\bullet})=
  \hom_{\cal{D}(\Db(A,B))}(M\bul,M^{\prime\bullet}).
$$
\dok
Since  $\overline{L}$ is an equivalence of categories,
\begin{gather*}
 \hom_{{\cal D}(Db(A,B))}(M\bul,M^{\prime\bullet})\\=
 \hom_{{\cal BAR}(\Db(A,b))}(\Barb(\Db(A,B),k,M\bul),
 \Barb(\Db(A,B),k,M^{\prime\bullet}))\\
 =H^0\left(\Hom_{\Db(A,B)}(\Barb(\Db(A,B),k,M\bul),
 \Barb(\Db(A,B),k,M^{\prime\bullet}))\right)_0.
\end{gather*}
Here $(*)_0$ denotes the zero grading component.

By the same arguments as in~\ref{bar0} it is proved that the
functor
\begin{gather*}
  X\bul\mapsto\Hom_{\Db(A,B)}(\Barb(\Db(A,B),k,M\bul),X\bul),\\
  \hom^n_{\Db(A,B)}(\Barb(\Db(A,B),k,M\bul),X\bul)=\prod_{q-p=n}
  \hom_{\Db(A,B)}(\Bar^p(\Db(A,B),k,M\bul),X^q)
\end{gather*}
maps exact DG-modules $X\bul\in\Db(A,B)\mod$
into exact complexes of vector spaces.

Thus
 \begin{multline*}
   \hom_{{\cal BAR}(\Db(A,B))}(\Barb(\Db(A,B),k,M\bul),
   \Barb(\Db(A,B),k,M^{\prime\bullet})) \\
   =\Hom_{{\cal BAR}(\Db(A,B))}(\Barb(\Db(A,B),k,M\bul),M^{\prime\bullet}).
   \qed
 \end{multline*}

Recall that the complex $K\bul=\Hom(\Barb(A,B,A),\underB)$
is both a {\em left\/} DG-module over $\Db(A,B)$ and a left
$A$-module. Denote the category of left DG-modules over
$\Db(A,B)$ satisfying the condition (D) by
$\dC(\Db(A,B)^{\opp})$.  Evidently
$K\bul\in\dC(\Db(A,B)^{\opp})$.

\sssn
\Lemma \label{inclusion}%
The natural morphism of algebras $A\map
\ext^0_{\Db(A,B)^{\opp}}(K\bul,K\bul)$ is an inclusion.
\ms
\dok
Consider the bar resolution $\Bar(\Db(A,B),k,K\bul)$.
By the previous Proposition
$$
  \ext^0_{\Db(A,B)^{\opp}}(K\bul,K\bul)=
  H^0\left(\Hom_{\Db(A,B)^{\opp}}(\Barb(\Db(A,B),k,K\bul),
  \Barb(\Db(A,B),k,K\bul))\right).
$$
Consider a functor
$$
  F:\ \dC(\Db(A,B))\map
  {\cal K}om(A\mod):\ X\bul\mapsto
  X\bul\ten_{\Db(A,B)}\Barb(\Db(A,B),k,K\bul)
$$
One checks directly that $F$ defines a functor on the
corresponding derived categories that is equivalent to
$\overline{D}^{\downarrow}_A$.

Suppose there exists $a\in A$ such that $a\mapsto 0$ in
$\ext^0_{\Db(A,B)^{\opp}}(K\bul,K\bul)$.  Then the
multiplication by $a$ defines an endomorphism of
$\Barb(\Db(A,B),k,K\bul)$ that is homotopic to zero, i.~e.
there exists a $\Db(A,B)$-linear map
$$
 h:\ \Barb(\Db(A,B),k,K\bul)\map\Bar^{\bullet-1}(\Db(A,B),k,K\bul)
$$
(not commuting with the differential) such that
$\operatorname{d}h+h\operatorname{d}=\mu_a$. Here $\mu_a$
denotes the endomorphism corresponding to $a\in A$.

Whence for every $X\bul\in\dC(\Db(A,B))$ the homomorphism of
complexes $F(X\bul)\map F(X\bul)$ given by the
multiplication by $a$ is homotopic to zero, in particular
$a\in A$ acts trivially on the cohomology spaces
$H\bul(F(X\bul))$.

The comparison of this fact with Theorem~\ref{equiv}
finishes the proof of the Lemma. \qed

\ssn
Now we are going to introduce an algebra $\oppA$ such that
$\oppA=B^{\opp}\ten N ^{\opp}$ as a vector space and
$\Db(A,B)^{\opp}\cong \Db(\oppA,B^{\opp})$.

Recall that the DG-algebra $\Db(A,B)$ is a tensor product of
the DG-subalgebra $\Db(N,k)$ and the subalgebra (not a
DG-subalgebra) $B^{\opp}$. Recall also that the algebra
$\Db(N,k)$ is isomorphic to the tensor algebra over the
graded vector space $\overline{N}^*$ and the differential in
it is generated by the map
$\overline{N}^*\map\overline{N}^*\ten\overline{N}^*$ dual to
the multiplication map in $N$ and extended to the whole
algebra $T(\overline{N}^*)$ by the Leibnitz rule.

To define a DG-algebra $C^{\vee}$ such that $C^{\vee}$ is a
tensor product of its DG-subalgebra $T(\overline{N}^*)$ with
a differential given by the map dual to the multiplication
in $N$ and a (not DG-) subalgebra $B^{\opp}$ one has to
specify the following additional data:
\begin{itemize}
  \item a linear map $B\ten
  \overline{N}^*\map\overline{N}^*\ten B$ generating the
  multiplication in $C^{\vee}$; \item a linear map $B\map
  \overline{N}^*\ten B$ providing a component of the
  differential in $C^{\vee}$,
\end{itemize}
satisfying certain
constraints (that provide the associativity constraint and
the Leibnitz rule in the DG-algebra $C^{\vee}$).

On the other hand to define an algebra $C$ such that $C$ is
a tensor product of two its subalgebras $B$ and $N$ as a
vector space and the conditions~\ref{setup} are satisfied
one has to specify the following data (additional to the
algebra structures on $B$ and $N$):
\begin{itemize}
  \item a linear map $\overline{N}\ten B\map B\oplus
  B\ten\overline{N}$ providing the multiplication in $C$
\end{itemize}
satisfying certain constraints (that provide
the associativity constraint in the algebra $C$).

\sssn
\Prop
The construction of the dual algebra $C\mapsto
C^{\vee}:=\Db(C,B)$ provides a one to one correspondence
between the two described types of data, i.~e. for every
DG-algebra $C^{\vee}$ such that
$C^{\vee}=T(\overline{N}^*)\ten B^{\opp}$ as a vector space
there exists an algebra $C=B\ten N$ as a vector space such
that the DG-algebras $C^{\vee}$ and $\Db(C,B)$ are
isomorphic.
\ms
\dok
Evidently the two described types of
data are dual  to each other.  One has to check directly
that the constraints in the first and the second case are
equivalent, i.~e. that the data of the first type defines a
DG-algebra (with the associativity constraint and the
Leibnitz rule satisfied) if and only if the dual data of the
second type defines an associative algebra.  \qed

\sssn
\Lemma
The DG-algebras $\Db(N,k)^\opp$ and $\Db(N^\opp,k)$ are
isomorphic to each other.  \qed

\sssn
\Cor
There exists an associative algebra $\oppA$ such that
$\oppA$ contains two subalgebras $B^{\opp}$ and $N^{\opp}$
satisfying  the conditions~\ref{setup} for $\oppA$, $B^\opp$
and $N^\opp$ such that the DG-algebra $\Db(A,B)^{\opp}\cong
\Db(\oppA,B^{\opp})$.  \qed

\sssn
\Cor
The functor $D_{\oppA}^\downarrow$ provides an equivalence
of triangulated categories \\
$\dcD(\oppA)\til{\map}\dcD(\Db(A,B)^\opp)$.  \qed

Now we describe the algebra $\oppA$ explicitly.

\ssn
Consider a left graded $N$-module $N^*:=\underset{n\ge
0}{\bigoplus} \hom_k(N_n,k)$.  The action of $N$ on the
space is defined as follows.
$$
 f:\ N\map k,\ n\in N,\text{
 then } (n\cdot f)(n'):=f(n'n).
$$
Consider also the left
$A$-module $S_A:=\Ind_N^A(N^*)=A\ten_NN^*$.  Evidently
$S_A\cong B\ten N^*$ as a $B$-module and by~\ref{setup}
$S_A\cong N^*\ten B$ as a $N$-module.

There is another left $A$-module with these two properties.
$S'_A:=\hom_B(A,B)$ with the left action of $A$ defined as
follows.
$$
 f:\ A\map B,\ a\in A,\text{ then }(a\cdot
 f)(a'):=f(a'a).
$$
\sssn
\Lemma
The $A$-modules $S_A$ and $S'_A$ are isomorphic.
\ms
\dok
Fix the decomposition $A=B\ten N$ provided by the
multiplication in $A$.  For $f\in N^*$ define $f_2:\ A\map
B$, $f_2(b\ten n):=f(n)b$.  Define the pairing
$$
 S_A\times A\map B,\ f\ten a\times a'\mapsto f_2(a'a).
$$
One checks directly the correctness of the definition and the
perfectness of the defined pairing.  \qed

Thus $S_A\cong \Coind_B^AB$. The functors of induction from
$N$ to $A$  and of coinduction from $B$ to $A$ provide the
natural inclusions of algebras
$$
 N^{\opp}\hookrightarrow
 \hom_A(S_A,S_A)\text{ and }B^\opp\hookrightarrow
 \hom_A(S_A,S_A).
$$
\sssn
\Lemma (i) $\hom_A(S_A,S_A)=\hom_k(N^*,B)$ as a vector
space;

(ii) the subalgebra in $\hom_A(S_A,S_A)$ generated by the
images of $N^{\opp}$ and $B^{\opp}$ is isomorphic to
$B^{\opp}\ten N^{\opp}$ as a vector space.  \qed

The latter algebra is denoted by $A^{\natural}$.

\sssn
\Theorem   \label{A=A}%
The algebras $\oppA$ and $A^{\natural}$ are isomorphic.
\ms
\dok
Consider the DG-module $\dD(S_A)$.

\sssn
\Lemma The right $\Db(A,B)$ DG-modules $\dD(S_A)$ and
$\Ind_{\Db(N,k)}^{\Db(A,B)}\underk$ are quasiisomorphic to
each other.  \qed

Consider also the left  DG-module over $\Db(A,B)^\opp$, that
is, the right DG $\Db(A,B)$-module $K_{\oppA}\bul$.

\sssn
\Lemma
$K_{\oppA}\bul$ is quasiisomorphic to
$\Ind_{\Db(N,k)}^{\Db(A,B)}\underk$.  \qed

By Theorem~\ref{equiv} the algebra $\hom_A(S_A,S_A)$ is
isomorphic to $\ext^0_{\Db(A,B)}(\dD(S_A),\dD(S_A))$. Hence
$A^{\natural}\hookrightarrow
\ext^0_{\Db(A,B)}\left(\Ind_{\Db(N.k)}^{\Db(A,B)}\underk,
\Ind_{\Db(N,k)}^{\Db(A,B)}\underk\right)$.

By Lemma~\ref{inclusion} $\oppA\hookrightarrow
\ext^0_{\Db(A,B)}\left(\Ind_{\Db(N.k)}^{\Db(A,B)}\underk,
\Ind_{\Db(N,k)}^{\Db(A,B)}\underk\right)$.

One can check directly that the images of the subalgebras
$B^\opp$ (resp.  $N^\opp$) under these two inclusions
coincide.  \qed

\sssn
\Rem
The previous Theorem provides an explicit way to recover the algebra
$C=B\ten N$ from $C^{\vee}=B^{\opp}\ten\Db(N,k)$. Namely consider the
{\em left\/} $C^{\vee}$ DG-module $\Ind_{\Db(N,k)}^{C^\vee}\underk$.
Since $\ext^0_{\Db(N,k)}(\underk,\underk)=N$ (this is a very easy
case of Theorem~\ref{equiv}), the algebra $N$ acts on
$\Ind_{\Db(N,k)}^{C^\vee}\underk$ by endomorphisms. Since
$\Db(N,k)$ is normal in $C^\vee$, i.~e. the left and the right
ideals in $C^\vee$ generated by its augmentation ideal coincide,
the algebra $B$ also acts on $\Ind_{\Db(N,k)}^{C^\vee}\underk$ by 
endomorphisms.

Then the subalgebra in
$\ext^0_{C^\vee}\left(\Ind_{\Db(N,k)}^{C^\vee}\underk,
\Ind_{\Db(N,k)}^{C^\vee}\underk\right)$ generated by the images of
$B$ and $N$ is isomorphic to $C$.

\sssn
\Lemma
Let $C_{(1)}\bul$ and $C_{(2)}\bul$ be two DG-algebras,
$$
  C_{(1)}\bul=\bigoplus_{p,q\in\Z}C_{(1)p}^q,\
  C_{(2)}\bul=\bigoplus_{p,q\in\Z}C_{(2)p}^q,
$$
and both $C_{(1)}\bul$ and $C_{(2)}\bul$ satisfy the condition (D).
Let $\varphi:\ C_{(1)}\bul\map C_{(2)}\bul$ be a quasiisomorphism of 
DG-algebras.
Then the restriction functor  provided by $\varphi$ is an
equivalence of triangulated categories
$$
  \dcD(C_{(2)}\bul)\til{\map}\dcD(C_{(2)}\bul).\qed
$$
In particular the algebra $C$ can be recovered from the
quasiisomorphism class of the DG-algebra $C^\vee$. \label{recover}%

\ssn
Consider a subcategory
${\cal C}^b(\Db(A,B))\hookrightarrow\Db(A,B)\mod$ that
consists of DG-modules satisfying both conditions (U) and
(D). Its localization by the class of quasiisomorphisms is
denoted by ${\cal D}^b(\Db(A,B))$.

\sssn
\Lemma Both inclusion functors
$$
  I^\uparrow:\ {\cal D}^b(\Db(A,B))\map \upcD(\Db(A,B))\text{
  and }I^\downarrow:\ {\cal D}^b(\Db(A,B))\map \dcD(\Db(A,B))
$$
are full and faithful.
\ms
\dok
Like the statement of Lemma~\ref{truncation}, the Lemma is an
easy consequence of the existence of truncation functors.
\qed

Consider a full triangulated subcategory
${\cal D}({\cal V}erma^*(A))$ in $\dcD(A)$ (finitely)
generated by all left $A$-modules of the form
 $\Coind_B^AL_0$, where $L_0$ is some finite dimensional
$B$-module, and a full triangulated subcategory
${\cal D}({\cal V}erma(A^{\sharp\opp}))$ in
$\upcD(A^{\sharp\opp})$ generated by all right
$\oppA$-modules of the form $\Ind_B^{A^{\sharp\opp}}L_0$,
where $L_0$ is some finite dimensional $B$-module.

\sssn
\Prop The functor $\overline{D}^{\uparrow}_{\oppA}\circ
D^\downarrow_A$ provides an equivalence of triangulated
categories ${\cal D}({\cal V}erma^*(A))$ and
${\cal D}({\cal V}erma(\oppA))$.
\ms
\dok
Note that for any
left $A$-module of the form $M=\Coind_B^AL_0$, where the
$B$-module $L_0$ is finite dimensional, $\dD_A(M)$ belongs
to ${\cal C}^b(\Db(A,B))$ and for any right $\oppA$-module of
the form $M'=\Ind_B^{A^{\sharp\opp}}L_0$, where the
$B$-module $L_0$ is finite dimensional, $\upD_{\oppA}(M')$
belongs to
${\cal C}^b(\Db(\oppA,B)^\opp)\til{\map}{\cal C}^b(\Db(A,B))$.

Now compare Theorem~\ref{equiv} with the previous Lemma.
\qed

\sssn
\Cor                         \label{vermamod}
The additive category of left  $A$-modules of the form
$\Coind_B^AL_0$ is equivalent to the additive category of
right $\oppA$-modules of the form
 $\Ind_B^{A^{\sharp\opp}}L_0$, where $L_0$ is finite
dimensional.  \qed

\subsection{Semiinfinite cohomology of associative algebras}
Here we give a definition of associative algebra
semiinfinite cohomology amd compare it with the one
presented in [Ar].

\sssn
\Def
Let $M\bul\in\upC(A^\opp)$,
$M^{\prime\bullet}\in\dC(\oppA)$. Then  set
$$
  \Exts_A(M^{\prime\bullet},M\bul):=
  \Ext_{\Db(A,B)^\opp}(\dD_{\oppA}(M^{\prime\bullet}),\upD_A(M\bul)).
$$
Note that by definition the semiinfinite $\ext$ functor
maps complexes exact by either of the variables to zero.

The definition of semiinfinite cohomology in [Ar] used the following standard
complex.
$$
  \stand(M^{\prime\bullet},M\bul):=
  \Hom_{\oppA}\left(\Barb(\oppA,B^\opp,M^{\prime\bullet}),
  \Barb(A^\opp,N^\opp,M\bul)\ten_AS_A\right).
$$
\sssn
\Theorem
Let $M\bul\in\upC(A^\opp)$,
$M^{\prime\bullet}\in\dC(\oppA)$. Then
$\Exts_A(M^{\prime\bullet},M\bul)=H\bul(\stand(M^{\prime\bullet},M\bul))$.
\ms
\dok
First we present several Lemmas.

\sssn
\Lemma \label{a}%
$S_{\oppA}\cong S_A$ as an $A-\oppA$ bimodule.
\qed

\sssn
\Lemma           \label{b}%
$\Barb(\Db(A,B)^\opp,\Db(N^\opp,k),\dD_{\oppA}(M^{\prime\bullet})\cong
\dD_{\oppA}(\Barb(\oppA,N^\opp,M^{\prime\bullet}))$ as a
left $\Db(A,B)$ DG-module.  \qed

\sssn
\Lemma                     \label{c}%
$K\bul_{\oppA}\cong \dD_A(S_A)$ as a $\Db(A,B)-\oppA$
DG-bimodule.  \qed

Let us construct an explicit isomorphism of complexes
$$
 \stand(M^{\prime\bullet},M\bul)\cong
 \Hom_{\Db(A,B)^\opp}\left(\Barb
 (\Db(A,B)^\opp,\Db(N^\opp,k),
 \dD_{\oppA}(M^{\prime\bullet})),\upD_A(M\bul)\right).
$$
Using Lemma~\ref{b} and recalling  that the functor
$\overline{D}_{\oppA}^{\downarrow}$ is right adjoint to
$\dD_{\oppA}$ we see that
$$
 RHS=
 \Hom_{\oppA}\left(\Barb(\oppA,N^\opp,M^{\prime\bullet}),
 K_{\oppA}\bul\ten_{\Db(A,B)}\upD_A(M\bul)\right).
$$
As before we
use the left $\oppA$-module structure on $K\bul_{\oppA}$.
Note also that as a graded module over $\Db(A,B)$ the
complex $\dD_A(S_A)$ is equal to $\Ind_B^{\Db(A,B)}S_A$.
Thus using Lemma~\ref{c} we obtain that
$$
  K_{\oppA}\bul\ten_{\Db(A,B)}\upD_A(M\bul)\cong
  \Barb(A^\opp,B^\opp,M\bul)\ten_AS_A
$$
as a graded $\oppA$-module (with a forgotten differential).
One can check directly that this isomorphism is an isomorphism of
complexes. Finally we have
$$
  RHS=\Hom_{\oppA}\left(\Barb(\oppA,N^\opp,M^{\prime\bullet}),
  \Barb(A^\opp,B^\opp,M\bul)\ten_AS_A\right)=LHS.
$$
Finally by Lemma~\ref{HHH}
$$
  \Exts_A(M^{\prime\bullet},M\bul)=H\bul\left(
  \Hom_{\Db(A,B)^\opp}(\Barb(\Db(A,B)^\opp,\Db(N,k)^\opp,
  \dD_{\oppA} (M^{\prime\bullet})), \upD_A(M\bul))\right).
  \qed
$$

\section{The algebra $\oppA$ for universal enveloping algebras}

Consider a graded Lie algebra
$\g=\underset{n\in\Z}{\bigoplus}\g_n,\ \dim\g_n<\infty$.
Suppose that $\g$ has a one dimensional central extension
$\hatg$, i.~e.  $\hatg=\g\oplus kC$ as a vector space and
the element $C$  of degree zero  is central in $\hatg$.

The bracket in $\hatg$ is given by the formula
$$
 [g_1,g_2]=[g_1,g_2]_{\g}+\omega(g_1,g_2)C,\ [g,C]=0,\
 g_1,g_2,g\in\g,
$$
where $\omega$ is some 2-cocycle of the
Lie algebra $\g$.

\ssn
Consider an associative algebra $A:=U(\hatg)/\{C-1\}$.
Here $U(\hatg)$
denotes the universal enveloping algebra of the Lie algebra
$\hatg$.

Then the grading on $\hatg$ induces  gradings both on its
universal enveloping algebra and on the algebra $A$ as the
relation $C=1$ is homogeneous.  Fix the natural triangular
decomposition: $A=B\ten N$ as a vector space, where
$$
  \b:=\underset{n\le 0}{\bigoplus}\g_n,\ B:=U(\b),\
  \n:=\underset{n>0} {\bigoplus}\g_n,\ N:=U(\n).
$$
Note that
the algebra $A$ with this triangular decomposition satisfies
the conditions~\ref{setup}.  Our main task in this section
is to describe the algebra $\oppA$ for the algebra $A$
explicitly.

\ssn
First recall the construction of the critical cocycle of the
Lie algebra $\g$.

\sssn
Let $V=\underset {n\in\Z}\bigoplus V_n$ be a graded vector
space such that $\dim V_n<\infty$ for $n>0$. Denote the
space $\underset{n\le 0}{\bigoplus}V_n$ (resp. the space
$\underset{n>0}{\bigoplus}V_n$) by $V_-$ (resp. by $V_+$).
Consider the Lie algebra $\gl(V)$ of linear transformations
of $V$ satisfying the condition:
$$
 a\in\gl(V),\
 a=(a_{ij}), \Longrightarrow \exists m\in\N:\ a_{ij}=0 \text{
 for }|i-j|>m.
$$
$\gl(V)$ has a well known central
extension given by the cocycle $\omega_0$:
$$
 \omega_0(a_1,a_2)=\tr(\pi\circ[a_1,a_2]-[\pi\circ
 a_1,\pi\circ a_2]).
$$
Here $\pi$ denotes the projection
$V\map V_+$ with the kernel $V_-$. Note that the trace is
well defined on maps $V_+\map V_+$.

\sssn
The adjoint representation of $\g$ provides the Lie algebra
morphism $\g\map\gl(\g)$. The inverse image of the 2-cocycle
$\omega_0$ is called the critical cocycle of the Lie algebra
$\g$ (see e.~g. [FFr]). We also denote it by $\omega_0$.

\ssn
Now we introduce an analogue of the DG-algebra $\Db(A,B)$.
As before denote the left $A$-module $A\ten_N\underk$ by
$\underB$.

\sssn
Consider a complex $\Db(\hatg,\b):=\Lambda(\n^*)\ten \underB$
with the differential defined as follows:
\begin{gather*}
  d:\ \Lambda^n(\n^*)\ten\underB\map\Lambda^{n+1}(\n^*)\ten\underB:\
  a_1,\ldots,a_{n+1}\in\n,\
  f:\ \Lambda^n(\n)\map\underB,\\
  df(a_1\wedge\ldots\wedge a_{n+1}):=
  \sum_{s=1}^{n+1}(-1)^s a_s\cdot f(a_1\wedge\ldots\wedge
  a_{s-1}\wedge a_{s+1} \wedge\ldots\wedge a_{n+1})\\
  +\sum_{s<t}(-1)^{s+t}f([a_s,a_t]\wedge a_1\wedge\ldots\wedge
  a_{s-1} \wedge a_{s+1}\wedge\ldots\wedge a_{t-1}\wedge
  a_{t+1}\wedge\ldots\wedge a_{n+1}).
\end{gather*}
We define
a  structure of the universal enveloping algebra of a graded
Lie superalgebra on $\Db(\hatg,\b)$.

Denote by $\underline{\n}$ the representation of the Lie
algebra $\b$ in the vector space $\n=\g/\b$. The dual
representation is denoted by $\underline{\n}^*$.

\sssn
Consider a graded Lie superalgebra $\a=\a^0\oplus \a^1$,
where $\a^0=\b$, $\a^1=\underline{\n}^*$, the bracket in
$\a$ is defined as follows: $\a^0\wedge\a^0\map\a^0$ is just
the bracket in $\b$, $\a^0\ten \a^1\map \a^1$ is given by
the representation of $\b$ in $\underline{\n}^*$.

Evidently $U(\a)\cong\Db(\hatg,\b)$ as a graded vector space.

\sssn
\Lemma
The differential in $\Db(\hatg,\b)$ satisfies the Leibnitz
rule.  \qed

\sssn
Consider a morphism of DG-algebras $\Db(A,B)\map\Db(\hatg,\b)$ defined
as follows. Consider first the canonical DG-algebra map
$\Db(U(\n),k)\map\Db(\n,0)$ provided by the inclusion
$\Lambda(\n)\hookrightarrow T(\n)$:
$$
  n_1\wedge\ldots\wedge n_s\mapsto\sum_{\sigma}\operatorname{sgn}
  (\sigma)n_{\sigma(n_1)}\ten\ldots\ten n_{\sigma(n_s)},
$$
where $\operatorname{sgn}(\sigma)$ denotes the sign of the transposition
$\sigma$. Note that the morphism commutes with differentials as it is induced
by the morphism of standard resolutions --- the bar resolution and the Lie one:
 \begin{gather*}
  U(\n)\ten\Lambda(\n)\ten\underk\map\Barb(U(\n),k,\underk),\\
  \Hom_{U(\n)}(\Barb(U(\n),k,\underk),\underk)\map
  \Hom_{U(\n)}(U(\n)\ten\Lambda(\n)\ten\underk,\underk).
 \end{gather*}
Now the induction functor provides the required morphism of DG-algebras
 \begin{multline*}
  \varphi:\ \Db(A,B)\cong\Hom_A(\Ind_N^A\Barb(N,k,\underk),\Ind_N^A\underk)\\
  \map  \Hom_A(\Ind_N^A(N\ten\Lambda(\n)\ten\underk),\Ind_N^A\underk)\cong
  \Db(\hatg,\b).
 \end{multline*}
\sssn
\Lemma
$\varphi$ is a quasiisomorphism of DG-algebras.
\qed

\Rem
Both the construction of $\Db(A,B)$ and the one of $\Db(\hatg,\b)$ are
particular cases of the general quadratic-linear Koszul duality
construction over a base ring due to L.~Positselsky [P]. From that point
of view the DG-algebra $\Db(A,B)$ is the dual object for the algebra
$A$ equipped with a filtration induced by the filtration $k\subset N$ on $N$,
and $\Db(\hatg,\b)$ is the dual object for the
algebra $A$ equipped with a filtration induced by the standard PBW 
filtration on
$N=U(\n)$.  The obvious morphism of filtered algebras leads to the 
morphism of
the dual objects $\varphi$.

\sssn
\Cor
(see~\ref{recover}) 
Let $\underB$ be the {\em left\/} $\Db(\hatg,\b)$ DG-module
$\Ind_{\Db(\n,0)}^{\Db(\hatg,\b)}\underk$. Then
$A$ coincides is isomorphic to the subalgebra in
$\ext^0_{\Db(\hatg,\b)}(\underB,\underB)$ generated by the
images of $B$ and $N$.
\qed

In particular the algebra $\oppA$ can be recovered not only from
$\Db(A,B)^\opp$ but also from $\Db(\hatg,\b)^\opp$. Now we calculate the
DG-algebra $\Db(\hatg,\b)^\opp$.

\ssn
Consider the antipode map in the Lie superalgebra $\a$:
$$
  \alpha:\ U(\a)\map U(\a)^\opp,\ \alpha(x)=-x,\ \alpha(y)=y
  \text{ for }x\in\a^0,\ y\in\a^1.
$$
Note that $\alpha$ need not necessarily commute  with the differential in
$U(\a)$. So to describe the algebra $\Db(\hatg,\b)^\opp$
it is sufficient to calculate relations between $\alpha$ and the differential.

Choose a homogeneous basis $\{x_m^{(i)}\}$ in $\g$, here $x_m^{(i)}\in
\g_i$. In particular the set $\{x_m^{(j)}\}_{j\le 0}$ forms a basis in $\b$,
the set $\{x_m^{(j)}\}_{j>0}$ forms a basis in $\n$. Let
$\{x_m^{(j)*}\}_{j>0}$ be the dual basis in $\n^*=\sum_{j>0}(\n_j)^*$.

\sssn
\Lemma

\qquad(i) For $x\in\Lambda(\n^*)$ the antipode map commutes with the
differential;

\qquad(ii) For $b\in\b$ the antipode map satisfies the relation
$$
  \alpha\cdot d(b)=d\cdot \alpha(b)-\sum_{n\in\Z,i>0}
  \omega_0(x_n^{(i)},b)x_n^{(i)*}.
$$
\dok
The first statement is proved by an obvious calculation. To prove the second
 one
note that
$$
    \alpha\cdot d(b)=\sum_{m\in\Z,j>0}\alpha(\pi_{\b}\circ[x_m^{(j)},b]\ten 
    x_m^{(j)*})
    =d\cdot\alpha(b)+\sum_{m\in\Z,j>0}(\pi_{\b}\circ[x_m^{(j)},b])
    (x_m^{(j)*}).
$$
The second summand is a linear function on $\n$. Let us calculate its value 
on a
base vector~$x_n^{(i)}$.
$$
  (\pi_{\b}\circ[x_m^{(j)},b])x_m^{(j)*}(x_n^{(i)})
  =x_m^{(j)*}([\pi_{\b}\circ[x_m^{(j)},b],x_n^{(i)}])=
  \tr(\pi_{\n}\circ[x_n^{(i)},\pi_{\b}\circ[b,*]]).
$$
The brackets here are the brackets in $A$, $\pi_{\n}$ and $\pi_{\b}$
denote the canonical projections onto $\n$ and $\b$ respectively.
Thus the second summand is equal to
$$
  \sum_{n\in\Z,i>0}\tr(\pi_{\n}\circ[x_n^{(i)},\pi_{\b}\circ[b,*]])x_n^{(i)*}
  =:\sum_{n\in\Z,i>0}\nu(x_n^{(i)},b)x_n^{(i)*}.
$$
On the other hand
 \begin{multline*}
  \nu(x_n^{(i)},b)=\tr(\pi_{\n}\circ\ad_{x_n^{(i)}}(\ad_b-\pi_{\n}
  \circ\ad_b))\\=
  -\tr(\pi_{\n}\circ([\ad_{x_n^{(i)}},\ad_b])-[\pi_{\n}\circ
  (\ad_{x_n^{(i)}}),\pi_{\n}\circ
  (\ad_b)])=-\omega_0(x_n^{(i)},b).\qed
 \end{multline*}
\sssn
\Cor
Let $\til{\g}$ be the central extension of the Lie algebra $\g$
defined with the help of the cocycle $-\omega_0-\omega$.
Then the algebra $\oppA$ is isomorphic to  $U(\til{\g})/\{C-1\}$.\qed

\sssn
\Rem
The proof of the fact that the algebra
$\operatorname{End}_{U(\g)}(S_{U(\g)})$ differs from $U(\g)$ by a Lie
algebra 2-cocycle is contained essentially in [V]. The method used there is
different from ours and is probably easier. Yet our proof follows from the
general picture that explains the very appearance of the cocycle.  In
particular, when changed a little, our proof still holds in the case of
affine quantum groups.

\subsection{Critical cocycle in the case of affine Lie algebras.}
To know
the bimodule structure on $S_A$ one has to describe not only the cohomology
class of the critical cocycle but the critical cocycle itself.  In this
section we present a direct calculation for the critical cocycle $\omega_0$ 
in
the case of affine Lie algebras. It turns out that the cocycle depends on
the type of the triangular decomposition of the affine Lie algebra.

\sssn
From now on we suppose that our base field is $\CC$.
Fix a Cartan matrix $(a_{ij})_{i,j=1}^r$ of the finite type and consider the
corresponding semisimple Lie algebra $\overline{\g}$.
Recall that the affine Lie algebra $\g$ for the semisimple Lie algebra
$\overline{\g}$ is defined as a central extension of a loop algebra ${\cal
L}\overline{\g}:=\overline{\g}\ten\CC$ $[t,t^{-1}]$.  Namely,
$\g:={\cal L}\overline{\g}\oplus\CC K$, and the bracket in $\g$ is defined as
follows:
$$
 [g_1\ten t^n,g_2\ten t^m]=[g_1,g_2]\ten
 t^{n+m}+\delta_{n+m,0}nB(g_1,g_2)K,
$$
where
$g_1,g_2\in\overline{\g},n,m\in \Z$, and $B(\cdot,\cdot)$ denotes the
Killing form of $\overline{\g}$.

Denote by $\overline{R},\overline{R}^+$ and
$\overline\Sigma=\{\alpha_1,\ldots,\alpha_r\}$ the root system of
$\overline{\g}$, the set of positive roots and the set of simple roots
respectively. Let $X$ be the corresponding weight lattice, $(\cdot|\cdot)$
denotes the canonical symmetric bilinear form on $X$.

Thus $\overline{\g}=\overline{\g}^+\oplus \h\oplus \overline{\g}^-$, where
$\h$ denotes the Cartan subelgebra in $\overline{\g}$, and
$$
\overline{\g}^+=\underset{\alpha\in\overline{R}^+}
{\bigoplus}\overline{\g}_{\alpha},\
\overline{\g}^-=\underset{\alpha\in\overline{R}^-}
{\bigoplus}\overline{\g}_{\alpha}.
$$
Consider the Chevalley generators of the Lie algebra $\overline{\g}$:
$
e_i\in\overline{\g}_{\alpha_i}, f_i\in\overline{\g}_{-\alpha_i},\
i=1,\ldots,r.
$

\sssn
It is well known that the affine Lie algebra $\g$ is a Kac-Moody Lie
algebra.  Denote by $R,R^+$ and $\Sigma$ the root system of
$\overline{\g}$, the set of positive roots and the set of simple roots
respectively. Then $\Sigma=\overline{\Sigma}\cup\{\alpha_0\}$. Thus
$\g=\n^+\oplus\h\oplus\CC K\oplus \n^-$, and
$$
\n^+=\underset{\alpha\in R^+}{\bigoplus}\g_{\alpha}=\overline{\g}^+\oplus
\overline{\g}\ten t\CC[t],\
\n^-=\underset{\alpha\in
R^-}{\bigoplus}\g_{\alpha}=\overline{\g}^-\oplus \overline{\g}\ten
t^{-1}\CC[t^{-1}].
$$
Recall the definition of the Chevalley generators for $\g$. Set
$e_{\alpha_i}:=e_i$, $e_{-\alpha_i}:=f_i$,
$i=1,\ldots,r$, where $e_i,f_i$
are the Chevalley generators for $\overline{\g}$,
$[e_i,f_i]=h_i$. Let $\alpha_{\operatorname{top}}$
be the highest root of the root system $\overline{R}$.
Choose a vector $h_{\operatorname{top}}$ in
$[\overline{\g}_{\alpha_{\operatorname{top}}},
\overline{\g}_{-\alpha_{\operatorname{top}}}]\subset\h$ such that
$\alpha_{\operatorname{top}}(h_{\operatorname{top}})=2$.
Fix
$e_{\operatorname{top}}\in\overline{\g}_{\alpha_{\operatorname{top}}}$
and
$f_{\operatorname{top}}\in\overline{\g}_{-\alpha_{\operatorname{top}}}$
such that
$[e_{\operatorname{top}},f_{\operatorname{top}}]=h_{\operatorname{top}}$.
Then set
$$
e_{\alpha_0}:=f_{\operatorname{top}}\ten t,\
e_{-\alpha_0}:=e_{\operatorname{top}}\ten t^{-1}.
$$
\sssn
\Lemma (see [K] (7.4.1))\label{kac}
$$[e_{\operatorname{top}},f_{\operatorname{top}}]=\sum_{i=1}^ra_i^{\vee}h_i,
$$
where $a_1^{\vee},\ldots,a_r^{\vee}$ are the marks for the dual Dynkin diagram.
\qed

Consider the following two natural gradings on $\g$ and the corresponding 
triangular
decompositions of $\g$.

\sssn
Let $\deg_{(1)}(g\ten t^n)=n$. The corresponding triangular decomposition
looks as follows:
$$
\g=\g_{>0}\oplus \g_{\le 0},\
\g_{>0}:=\overline{\g}\ten t\CC[t],\ \g_{\le 0}:=\overline{\g}\ten\CC
[t^{-1}].
$$
\sssn
Let $\deg_{(2)}$ be the grading on the Lie algebra $\g$ obtained fron the
root decomposition
of $\g$ by putting
$$
\deg_{(2)}(g_{\alpha})=1,\text{ where }g_{\alpha}\in\g_{\alpha},
\text{ for any }
\alpha\in\Sigma.
$$
Then the corresponding triangular decomposition is as follows:
$$
\g=\n^+\oplus\b^-,\
\b^-:=\n^-\oplus\h\oplus\CC K.
$$
\sssn
\Rem More generally every parabolic subalgebra $\p$ in $\overline{\g}$ with
nilpotent radical $\gothu\subset\overline{\g}^+$
provides a triangular decomposition of $\g$ such that the positive Lie 
subalgebra
is equal to $\gothu\oplus\overline{\g}\ten t\CC[t]$. Thus here we consider
the cases of $\p$ equal to $\overline{\g}$ and to the Borel subalgebra in 
$\overline{\g}$.

Denote the critical cocycle that corresponds to the first (resp. the second)
triangular decomposition by $\omega_0^{(1)}$ (resp. by $\omega_0^{(2)}$).

\sssn
\Lemma
Let $\omega$ be a 2-cocycle of the Lie algebra $\g$ that preserves the
root grading. Then $\omega$
is completely defined by the set of its values
$\{\omega(e_{\alpha},e_{-\alpha})|\alpha\in\Sigma\}$.
\ms
\dok
By definition of a 2-cocycle for any $g_1,g_2,g_3\in\g$
$$
\omega([g_1,g_2],g_3)=\omega(g_1,[g_2,g_3]).
$$
Moreover the cocycle $\omega$ is nontrivial on a pair of homogeneous elements
$g_1,g_2\in\g$ only if $\g_1\in\g_{\alpha},g_2\in\g_{-\alpha},\alpha\in R$
or both $g_1$ and $g_2$ belong to $\h\oplus\CC K$.
By definition of the affine algebra
for any $\alpha\in R$ $\dim\g_{\alpha}=1$, and
for any $\alpha\in R^+$
there exist  $\beta_1,\ldots,\beta_s\in\Sigma$ such that
$$
e_{\alpha}:=[e_{\beta_1},[e_{\beta_2},\ldots,
[e_{\beta_{s-1}},e_{\beta_s}]\ldots]
\ne0\in\g_{\alpha}
$$
(resp.
for any $\alpha\in R^-$
there exist  $\beta_1,\ldots,\beta_s\in\Sigma$ such that
$$
e_{-\alpha}:=[e_{-\beta_1},[e_{-\beta_2},\ldots,[e_{-\beta_{s-1}},
e_{-\beta_s}]\ldots]
\ne0\in\g_{\alpha}).
$$
The vector space $\h$ has a base consisting of $h_i:=[e_{\alpha_i},
e_{-\alpha_i}],
\alpha_i\in\overline{\Sigma}$. The space $\CC K$ is generated by the vector
$[h_1\ten t,h_1\ten t^{-1}]$. Note that  $h_1\ten t\in\n^+,\ h_1\ten
t^{-1}\in\n^-$.  Now the statement of the Lemma is proved by induction by
the length of the expression of the homogeneous elebment $g_1$ via the
generators of $\g$.  \qed

Consider the 1-cochain $\rho$ of $\g$ defined as follows:
$$
\rho(g_{\alpha})=0,\ g_{\alpha}\in\g_{\alpha},\ \rho(K)=0,\ \rho(h_i)=1,\
i=1,\ldots,r.
$$
Then for every $\alpha_i\in\overline{\Sigma}\  d\rho(e_{\alpha_i},
e_{-\alpha_i})
=1$.

Recall that the dual Coxeter number  $h^{\vee}$ of the root system
$\overline{R}$ is equal to $a_1^{\vee}+\ldots+a_r^{\vee}+1$.

\sssn
\Lemma
$d\rho(e_{\alpha_0},e_{-\alpha_0})=h^{\vee}-1$.
\ms
\dok
Follows immediately from~\ref{kac}.\qed

Denote by $\nu$ the canonical 2-cocycle of $\g$
$$
(g_1\ten t^k,g_2\ten t^m)\mapsto n\delta_{k+m,0}(g_1,g_2), g_1,g_2\in
\overline{\g}, k,m\in\Z,
$$
where $(\cdot,\cdot)$ denotes the normalized invariant bilinear
form on $\overline{\g}$:
$$
(\cdot,\cdot):=\frac 1 {2h^{\vee}} B(\cdot,\cdot).
$$
\sssn
\Prop

\qquad(i) $\omega_0^{(1)}=2h^{\vee}\nu$

\qquad(ii) $\omega_0^{(2)}=2h^{\vee}\nu+2d\rho$.

\dok
(i) Since the cocycle $\omega_0^{(1)}$ preserves the first grading,
for every $g_1,g_2\in \overline{\g}$
$$
\omega_0^{(1)}(g_1\ten t,g_2\ten t^{-1}) =
\tr_{\overline{\g}\ten t} \left(\operatorname{ad}(g_1\ten t)\operatorname{ad}
(g_2\ten t_{-1})\right)=B(g_1,g_2).
$$
(ii) To calculate $\omega_0^{(2)}(e_{\alpha_i},e_{-\alpha_i})$, 
$i=1,\ldots,r$
note that as the cocycle preserves the root grading, the only nontrivial
summand
in the root sum will be
$$
\langle e_{\alpha_i}^*,[e_{\alpha_i},[e_{-\alpha_i},e_{\alpha_i}]]\rangle.
$$
To calculate $\omega_0^{(2)}(e_{\alpha_0},e_{-\alpha_0})$ note that the 
only
root vector $e$ in $\overline{\g}\ten t$ such that 
$[e_{-\alpha_0},e]\in\b^-$
is $e_{\alpha_0}$.
\qed

\sssn
Fix a level $k$ and a character $\lambda$ of $\h$ (resp a finite dimensional
representation of $\overline{\g}$ $L(\mu)$). Recall that the module
$M(\lambda):=\Ind_{U_(\n^+\oplus\h)}^{U_k(\g)}{\CC}(\lambda)$
is called the Verma module,
and the module $W(\mu):=\Ind_{U(\g_{\ge 0})}^{U_k(\g)}L(\mu)$
is called the Weyl
module over $\g$ on the level $k$.

\Rem
Comparing the previous  calculation with~\ref{vermamod} we obtain in 
particular the
following statement:

\qquad(i) The additive categories generated by Verma modules over $\g$ on 
levels
$k$ and $-2h^{\vee}-k$ are

\qquad antiequivalent;

\qquad(ii) The additive categories generated by Weyl modules over $\g$ on
levels $k$ and $-2h^{{\vee}-k}$ are

 \qquad antiequivalent.

\section*{References}

$\text{[Ar]}$ S.Arkhipov. {\it Semiinfinite cohomology of quantum
groups.} Preprint q-alg/9601026 (1996), 1-24.\\
$\text{[BGG]}$ J.N.Bernshtein, I.M.Gelgand, S.I.Gelfand. {\it
Algebraic bundles over \hbox{\bf P}$^n$  and problems of linear algebra.}
Functional Analysis and Appl. \hbox{\bf 12} (1978) no.~3, pp.~66-67.\\
$\text{[BGS]}$ A. Beilinson, V. Ginzburg, V. Schechtman.
{\it Koszul duality.} Journal of geometry and physics \hbox{\bf 5} (1988),
317-350.\\
$\text{[F]}$ B.Feigin.
{\it Semi-infinite cohomology of Kac-Moody and Virasoro Lie algebras.} Usp.
Mat.  Nauk \hbox{\bf 33},no.2, 195-196, (1984),  (in Russian).\\
$\text{[FFr]}$ B. Feigin, E. Frenkel. {\em Semi-infinite Weil complex and the
Virasoro algebra.} Comm. Math. Phys. \hbox{bf 137} (1991), 617-639.\\
$\text{[GeM]}$
S.I.Gelfand, Yu.I.Manin.  {\it Methods of homological algebra.} M:  Nauka,
 (1988).  (In Russian)\\
$\text{[K]}$ V. Kac. {\it Infinite dimensional Lie algebras.} Birkh\"auser,
Boston (1984).\\
$\text{[P]}$ L. Positselsky. {\it Quadratic-linear Kozsul duality over a
base ring}, unpublished.\\
$\text{[V]}$
A.Voronov. {\it Semi-infinite
homological algebra.} Invent. Math. {\bf 113}, (1993),  103--146.

\end{document}